\documentclass{amsart}

\usepackage{amsmath,amssymb}
\usepackage[all]{xy}

\newtheorem{theorem}{Theorem}
\newcommand{\bt}{\begin{theorem}}
\newcommand{\et}{\end{theorem}}

\newtheorem*{theoremNN}{Theorem}
\newcommand{\btNN}{\begin{theoremNN}}
\newcommand{\etNN}{\end{theoremNN}}

\newtheorem{lemma}{Lemma}
\newcommand{\bl}{\begin{lemma}}
\newcommand{\el}{\end{lemma}}

\newtheorem{corollary}{Corollary}
\newcommand{\bc}{\begin{corollary}}
\newcommand{\ec}{\end{corollary}}

\newtheorem{definition}{Definition}
\newcommand{\bdf}{\begin{definition}}
\newcommand{\edf}{\end{definition}}

\newtheorem{conjecture}{Conjecture}
\newcommand{\bconj}{\begin{conjecture}}
\newcommand{\econj}{\end{conjecture}}

\newtheorem*{conjectureNN}{Conjecture}
\newcommand{\bconjNN}{\begin{conjectureNN}}
\newcommand{\econjNN}{\end{conjectureNN}}

\newtheorem{example}{Example}
\newcommand{\bex}{\begin{example}}
\newcommand{\eex}{\end{example}}

\newtheorem{problem}{Problem}
\newcommand{\bprob}{\begin{problem}}
\newcommand{\eprob}{\end{problem}}

\newtheorem{oproblem}{Open Problem}
\newcommand{\boprob}{\begin{oproblem}}
\newcommand{\eoprob}{\end{oproblem}}

\newcommand{\beq}{\begin{equation}}
\newcommand{\eeq}{\end{equation}}
\newcommand{\benum}{\begin{enumerate}}
\newcommand{\eenum}{\end{enumerate}}

\newcommand{\N}{\ensuremath{ \mathbf N }}

\newcommand{\Z}{\ensuremath{\mathbf Z}}

\newcommand{\mcu}{\ensuremath{ \mathcal U}}
\newcommand{\mcv}{\ensuremath{ \mathcal V}}

\newcommand{\bq}{\begin{eqnarray*}}
\newcommand{\eq}{\end{eqnarray*}}
\newcommand{\be}{\begin{eqnarray}}
\newcommand{\ee}{\end{eqnarray}}
\newcommand{\ba}{\begin{array}}
\newcommand{\ea}{\end{array}}
\newcommand{\bfr}{\begin{flushright}}
\newcommand{\efr}{\end{flushright}}

\newcommand{\bmat}{\left(\begin{matrix}}
\newcommand{\emat}{\end{matrix}\right)}

\DeclareMathOperator{\qqand}{\qquad\text{and}\qquad}

\title[A partition theorem of Sylvester]{Trapezoidal numbers, divisor functions, 
and a partition theorem of Sylvester}
\author{Melvyn B. Nathanson}
\address{Department of Mathematics, Lehman College (CUNY), Bronx, NY 10468}
\email{melvyn.nathanson@lehman.cuny.edu}


\begin{document}
\maketitle

\begin{center} To Krishnaswami Alladi on his 60th birthday \end{center}

\begin{abstract}
A \emph{partition}  of  a positive integer $n$ is a representation of $n$ 
as a sum of a finite number of positive integers (called \emph{parts}).
A \emph{trapezoidal number} is a positive integer that has a partition whose parts 
are a decreasing sequence of consecutive integers, or, more generally, whose parts form 
a finite arithmetic progression.  This paper reviews the relation 
between trapezoidal numbers, partitions, and the set of divisors of a positive integer.  
There is also a complete proof of a theorem of Sylvester 
that produces a stratification of the partitions of an integer into odd parts 
and partitions into disjoint trapezoids.
\end{abstract}

\section{Partition theorems of Euler and Sylvester}
Let \N, $\N_0$, and \Z\ denote, respectively, the sets of positive integers, 
nonnegative integers, and integers.  
A \emph{partition} of  a positive integer $n$ is a representation of $n$ 
as a sum of a finite number of positive integers (called \emph{parts}), 
written in decreasing order.  
The usual left-justified Ferrers diagram of the partition 
\[
n = a_1 + a_2 + \cdots + a_k 
\]
with 
\[
a_1 \geq a_2 \geq \cdots \geq a_k \geq 1
\]
consists of $k$ rows of dots, 
with $a_i$ dots on row $i$.
For example, the Ferrers diagram of the partition 
\[
57 = 11 + 11 + 11 + 9 + 5 + 5 + 5
\]
is 
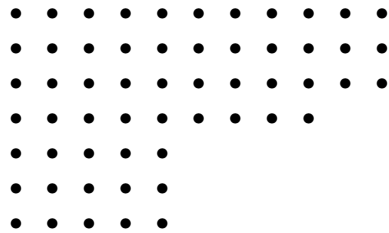
\[
\xymatrix@=0.1cm{
\bullet & \bullet & \bullet & \bullet & \bullet & \bullet & \bullet & \bullet &\bullet & \bullet & \bullet \\
\bullet & \bullet & \bullet & \bullet & \bullet & \bullet & \bullet & \bullet &\bullet & \bullet & \bullet \\
\bullet & \bullet & \bullet & \bullet & \bullet & \bullet & \bullet & \bullet &\bullet & \bullet & \bullet \\
\bullet & \bullet & \bullet & \bullet & \bullet & \bullet & \bullet & \bullet &\bullet  \\
\bullet & \bullet & \bullet & \bullet & \bullet \\
\bullet & \bullet & \bullet & \bullet & \bullet \\
\bullet & \bullet & \bullet & \bullet & \bullet 
}
\]

Perhaps the best known result about  partitions is the following theorem of Euler.  

\bt[Euler]          \label{trapezoid:theorem:Euler-odd-distinct}
The number of partitions of $n$ into odd parts equals 
the number of partitions of $n$ into distinct parts.
\et

\begin{proof}
Let $p_{odd}(n)$ denote the number of partitions of $n$ into odd parts, 
and let $p_{dis}(n)$ denote the number of partitions into distinct parts.  
A deceptively simple proof uses formal power series:
\begin{align*}
\sum_{n=0}^{\infty} p_{odd}(n) q^n 
& = \prod_{n=1}^{\infty} \frac{1}{1-q^{2n-1} } 
= \prod_{n=1}^{\infty} \frac{1-q^{2n} }{(1-q^{2n-1} )(1-q^{2n} )} \\
& = \prod_{n=1}^{\infty} \frac{1-q^{2n} }{1-q^n } 
= \prod_{n=1}^{\infty} \frac{ (1-q^n)(1+q^n)}{1-q^n } \\
& = \prod_{n=1}^{\infty} (1+q^n) 
= \sum_{n=0}^{\infty} p_{dis}(n) q^n.
\end{align*}
This argument is valid only after one understands infinite products, inversion, 
and composition of formal power series.  
\end{proof}

Every positive integer $n$ has a unique $g$-adic representation 
in the form $n = \sum_{i=0}^{\infty} \varepsilon_i g^i$, where 
$\varepsilon_i \in \{0, 1,\ldots, g-1\}$ for $i \in \N_0$ and $\varepsilon_i = 0$ 
for all sufficiently large $i$.  
Glaisher~\cite{glai83} generalized Euler's theorem by using the uniqueness 
of the $g$-adic representation.  
Theorem~\ref{trapezoid:theorem:Euler-odd-distinct} is the special case $g=2$.

\bt[Glaisher]          \label{trapezoid:theorem:Glaisher-odd-distinct}
Let $g \geq 2$.  The number of partitions of $n$  into parts not divisible by $g$ 
equals the number of partitions of $n$ such that every part occurs less than 
$g$ times.  
\et

\begin{proof}
Every positive integer $a$ can be written uniquely in the form $a = g^v s$, 
where $s$ is not divisible by $g$.  
Sylvester calls $s$ the \emph{nucleus} of $a$.
A partition of $n$ in which every part occurs at most $g-1$ times 
can be written uniquely in the form 
\beq          \label{trapezoid:odd1}
n = \varepsilon_1 a_1 + \cdots + \varepsilon_k a_k
\eeq
where the parts $a_1,\ldots, a_k$ are pairwise distinct  
and $\varepsilon_i \in \{1,\ldots, g-1\}$ for $i=1,\ldots, k$.  
Let
\[
a_i = g^{v_i} s_i = \underbrace{s_i + \cdots + s_i}_{\text{$g^{v_i}$ summands}}
\]
where $s_i$ is the nucleus of $a_i$.  
The nuclei $s_1, \ldots, s_k$ are not necessarily distinct.  
Let $S = \{s_1,\ldots, s_k\}$.  For each $s \in S$, let  
\[
\delta(s) =  \sum_{\substack{i\in \{1,\dots, k\}\\ s_i = s  }} \varepsilon_i g^{v_i}.
\]
Then 
\begin{align*}
n & = \varepsilon_1 a_1 + \cdots + \varepsilon_k a_k \\
& = \varepsilon_1 g^{v_1} s_1 + \cdots + \varepsilon_k g^{v_k} s_k \\
& = \underbrace{s_1 + \cdots + s_1}_{\text{$\varepsilon_1g^{v_1}$ summands}} 
+ \cdots + 
\underbrace{s_k + \cdots + s_k}_{\text{$\varepsilon_k g^{v_k}$ summands}} \\
& = \sum_{s\in S} \left(  \sum_{\substack{i\in \{1,\dots, k\}\\ s_i = s  }} 
\varepsilon_i g^{v_i}\right) s \\
& =  \sum_{s\in S} \delta(s) s.
\end{align*}
Thus, from the partition~\eqref{trapezoid:odd1}  of $n$ into parts occurring
less than $g$ times we have constructed a partition 
of $n$ as a sum of  integers not divisible by $g$. 

Conversely, let $n =  \sum_{s\in S}\delta(s) s $ be a partition of $n$ with 
parts in a set $S$ of integers not divisible by $g$, 
and where each $s \in S$ has multiplicity $\delta(s)$.  
Consider the $g$-adic representation 
\[
\delta(s) = \sum_{i\in I_s} \varepsilon_i g^i
\]
where $\varepsilon_i \in \{1,\ldots, g -1\}$.  
If $( i_1, s_1) \neq (i_2, s_2)$, then 
$g^{i_1}s_1 \neq g^{i_2}s_2$
and so 
\[
n =  \sum_{s\in S }\delta(s) s 
=  \sum_{s\in S } \sum_{i\in I_s} \varepsilon_i g^i  s 
\]
is a partition of $n$ into distinct parts $g^is$ with multiplicities at most $g-1$.
These two partition transformations are inverse maps, and
establish a one-to-one correspondence between partitions into parts 
not divisible by $g$ and parts occurring with multiplicities less than $g$.  
\end{proof}

Sylvester~\cite[sections 45--46]{sylv82} discovered and proved a different, very beautiful, 
and insufficiently known generalization 
of Euler's theorem.  We prove this theorem in 
Section~\ref{trapezoid:section:Sylvester}.

\section{Trapezoidal numbers}

For integers $k\in \N$,  $t\in \N_0$, and $a \in \Z$, 
the finite arithmetic progression with length $k$, difference $t$, 
and first term $a$ is the set
\beq     \label{trapezoid:ap}
\{a,a+t,a+2t,\ldots, a+(k-1)t \}.
\eeq
The sum of this arithmetic progression is 
\beq     \label{trapezoid:ska}
s_{k,t}(a) = \sum_{i=0}^{k-1} (a+it) = ka + \frac{k(k-1)t}{2}.
\eeq
The integer $a$ is the smallest element of the set~\eqref{trapezoid:ap} because $t \geq 0$.

Let $t  \in \N_0$.  A positive integer $n$ is a \emph{$k$-trapezoid with difference $t$}
if it is the sum 
of a finite arithmetic progression of  integers of length $k$ and difference $t$, 
that is, if it can can be represented in the form~\eqref{trapezoid:ska} 
for integers $k\in \N$, $t\in \N_0$, and $a \in \Z$.
A \emph{trapezoid with difference $t$} is a 
$k$-trapezoid with difference $t$ for some $k \in \N$.  
A \emph{$k$-trapezoid} is a $k$-trapezoid with difference $1$.  
For example, every odd integer is a 2-trapezoid, because $2n-1 = (n-1) + n$.  
A  \emph{trapezoid} is an integer that is a $k$-trapezoid for some $k$, that is, 
an integer that can be represented as the sum of a strictly decreasing sequence of consecutive 
integers.

A $k$-trapezoid with difference $t$ is \emph{positive} if $a \geq 1$ and 
\emph{nonpositive} if $a \leq 0$.  
If $a$ is positive, then the Ferrers diagram of this partition of $n$ 
has a trapezoidal shape.  
For example, $32 = 11 + 9 + 7 + 5$ is a positive 4-trapezoid with difference 2.  
Its Ferrers diagram is
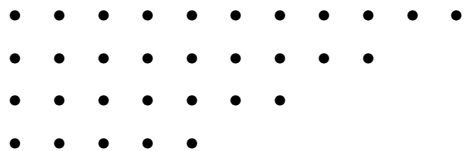
\[
\xymatrix@=0.2cm{
\bullet & \bullet & \bullet & \bullet & \bullet & \bullet & \bullet & \bullet & \bullet & \bullet & \bullet \\
\bullet & \bullet & \bullet & \bullet & \bullet & \bullet & \bullet & \bullet & \bullet \\
\bullet & \bullet & \bullet & \bullet & \bullet & \bullet & \bullet \\
\bullet & \bullet & \bullet & \bullet  & \bullet \\
}
\]

Every positive integer $n$ has a trivial positive trapezoidal representation 
with length 1 and difference 1, namely, $n = n$.  
Sylvester~\cite{sylv83} and Mason~\cite{maso12} proved 
that a positive integer $n$ is a $k$-trapezoid for some $k \geq 2$ 
if and only if $n$ is not a power of 2, and that  
the number of positive trapezoidal representations 
of $n$ is exactly the number of odd positive divisors of $n$.  
Bush~\cite{bush30} extended this result to 
trapezoidal representations with difference $t$.   
We prove their theorems below.

 In Section~\ref{trapezoid:section:Sylvester}  
we show how a special case of a partition theorem of Sylvester establishes 
another bijection between the number of trapezoidal representations of $n$ 
and the number of positive odd divisors of $n$.

For every positive integer $n$, let $\Phi_t(n)$ denote the number of representations 
of $n$ as a trapezoid with difference $t$, and let  $\Phi_t^+(n)$ denote the number 
of representations of $n$ as a positive trapezoid with difference $t$.
Thus,
\begin{align*}
\Phi_t(n) & = \left| \left\{  (k,a) \in \N \times \Z: s_{k,t}(a) = n  \right\} \right| \\ 
\Phi_t^+(n) & = \left| \left\{  (k,a) \in \N \times \N: s_{k,t}(a) = n  \right\} \right| 
\end{align*}
For $t = 1$, these functions count partitions into consecutive integers.  

Let $d(n)$ denote the number of positive divisors of $n$, 
and let $d_1(n)$ denote the number of odd positive divisors of $n$.
Let $d(n,\theta)$ denote the number of positive divisors $d$ of $n$ 
such that $d < \theta$.  If $n/2 < k \leq n$, then 
\[
d(n,k) = d(n,n) = d(n)-1.
\]
Let $[x]$ denote the integer part of the real number $x$.

\bl
Let $t$ and $n$ be  positive integers.
For every positive integer $k$, there is at most one representation of $n$ 
as a sum of a $k$-term arithmetic progression of integers with difference $t$.  
\el

\begin{proof}
This is true because the function 
$s_{k,t}(a)$ defined by~\eqref{trapezoid:ska} is a strictly increasing 
function of $a$.
\end{proof}

\bt                  \label{trapezoid:theorem:t-even}
Let $t$ be an even positive integer.  
For every positive integer $n$, 
\beq   \label{trapezoid:even-ell}
\Phi_t(n) = d(n)
\eeq
and 
\beq   \label{trapezoid:even+ell}
\Phi_t^+(n) = d(n, \theta)
\eeq
where
\[
\theta = \frac{1}{2} + \sqrt{\frac{2n}{t}+\frac{1}{4}}.
\]
\et

\begin{proof}
For every positive divisor $k$ of $n$, 
\[
a_{k,t}(n) = \frac{n}{k} - \frac{(k-1)t}{2}
\]
is an integer and  
\[ 
s_{k,t}\left( a_{k,t}(n) \right) = \sum_{i=0}^{k-1} \left(  \frac{n}{k} - \frac{(k-1)t}{2} +it \right) = n. 
\]
Moreover, if $k$ and $d$ are distinct positive divisors of $n$, 
then $a_{k,t}(n) \neq a_{d,t}(n)$.
Thus, $d(n) \leq \Phi_t(n) $. 

Conversely, if $n$ is the sum of a $k$-term arithmetic progression 
with even difference $t$ and first term $a$, then 
\[
n = s_{k,t}(a)   = k \left( a+\frac{(k-1)t}{2} \right)
\]
and so  $k$ is a positive divisor of $n$ and  $a = a_{k,t}(n).$
Thus, $\Phi_t(n) \leq d(n)$, and so there is a one-to-one correspondence 
between the positive divisors of $n$ 
and  representations of $n$ as a sum of a finite arithmetic progression 
with difference $t$.
This proves~\eqref{trapezoid:even-ell}.

Let $n = \sum_{i=0}^{k-1} (a+it)$.  The first term $a = a_{k,t}(n)$ is positive if and only if 
\[
\frac{n}{k} > \frac{(k-1)t}{2}
\]
or, equivalently,
\[
k < \frac{1}{2} + \sqrt{\frac{2n}{t}+\frac{1}{4}}.
\]
This proves~\eqref{trapezoid:even+ell}.
\end{proof}

\bl                  \label{trapezoid:lemma:t-odd}
Let $t$ be an odd positive integer.  
Let $n$ be a positive integer, and let $s_{k,t}(a) = n$ 
for some integer $a$ and some positive integer $k$.  
If $k$ is odd, then $k$ is an odd positive divisor of $n$.
If $k$ is even, then $2n/k$ is an odd positive divisor of $n$.
\el

\begin{proof}
If $k$ is odd, then $(k-1)/2$ is an integer and the identity 
\[
n = s_{k,t}(a) = ka + \frac{k(k-1)t}{2} = k\left(a + \frac{(k-1)t}{2} \right)
\]
implies that $k$ is a positive divisor of $n$.  

If $k$ is even, then $d = 2a + (k-1)t$ is odd and the identity 
\[
n =  \frac{k}{2}(2a + (k-1)t) 
\]
implies that $2n/k = 2a + (k-1)t$ is an odd positive divisor of $n$.  
This completes the proof.  
\end{proof}

\bt                     \label{trapezoid:theorem:t-odd}
Let $t$ be an odd positive integer.  
For every odd positive  divisor $k$ of $n$,  
there is exactly one representation of
$n$ as a sum of a $k$-term arithmetic progression of integers with difference $t$,  
and there is exactly one representation of $n$ as a sum 
of a $(2n/k)$-term arithmetic progression of integers with difference $t$.

The number of representations of $n$ as a $t$-trapezoid is 
\[
\Phi_t(n) = 2 d_1(n).
\]
\et

\begin{proof}
Let $k$ be a odd positive  divisor of $n$, and let $n = kq$.  
If $k = 2e+1$, then 
\beq   \label{trapezoid:odd}
n = \sum_{i=-e}^e (q+it) 
\eeq
is a representation of $n$ as a sum of an arithmetic progression 
with difference $t$, length $k$, and first term 
\beq   \label{trapezoid:a}
a_{k,t}(n) = q-et = \frac{n}{k} - \frac{(k-1)t}{2}.
\eeq

Let  
%\[ \ell = \frac{2n}{k} =2q \] and
\beq   \label{trapezoid:b}
b_{k,t}(n) = \frac{n}{2q}- \frac{(2q -1)t}{2} = \frac{k + t}{2} - \frac{nt}{k}.
% = \frac{k + t}{2}  - \frac{nt}{k} 
\eeq
Then 
%$\ell$ is an even positive integer, 
$b_{k,t}(n)$ is an integer, and 
\beq   \label{trapezoid:even}
n = \sum_{i=0}^{2q -1} (b_{k,t}(n)+it) 
\eeq
is a representation of $n$ as a sum of an arithmetic progression 
with difference $t$,  
length $2q = 2n/k$, and first term $b_{k,t}(n)$.
Applying Lemma~\ref{trapezoid:lemma:t-odd}, we see that there 
is a one-to-one correspondence between the odd positive divisors of $n$ 
and the representations of $n$ as a sum of an arithmetic progression 
with difference $t$ and odd length, and there is also 
a one-to-one correspondence between the odd positive divisors of $n$ 
and the representations of $n$ as a sum of an arithmetic progression 
with difference $t$ and even length.
This completes the proof.  
\end{proof}

For example, the only odd positive divisor of 1 is 1, 
and so $\Phi_t(1) = 2d_1(1) = 2$.   
The two representations of 1 as a sum 
of a finite arithmetic progression with odd difference $t$ are $1 = 1$ and 
\[
1 = \left( \frac{1-t}{2} \right) +  \left(  \frac{1+t}{2} \right).  
\]

The only odd positive divisor of 2 is 1, 
and so $\Phi_t(2) = 2d_1(1) = 2$.   
The two representations of 2 as a sum 
of a finite arithmetic progression with odd difference $t$ are $1 = 1$ and 
\[
2 = \left(\frac{1-3t}{2}\right)  + \left(\frac{1-t}{2}\right)  + \left(\frac{1+t}{2}\right)  + \left(\frac{1+3t}{2}\right).
\]

The trapezoidal representations with odd difference $t$ of an odd prime $p$ are 
\begin{align*}
p & = \frac{p-t}{2} + \frac{p+t}{2}  \\
& = \sum_{i=0}^{p-1} \left( 1 +    \frac{(2i - p + 1)t }{2}  \right) \\
& = \sum_{i=0}^{2p-1} \frac{ 1 + (2i - 2p+1)t }{2}.
\end{align*}
Thus, the four trapezoidal representations with difference 3 of the prime 5 are
\begin{align*}
5 & = 1 + 4 \\
& = (-5) + (-2) + 1 + 4 + 7 \\
& = (-13) + (-10) + (-7) + (-4) + (-1) + 2 + 5 + 8 + 11 + 14.
\end{align*}

\bt           \label{trapezoid:theorem:PowerOf2}
For every positive integer $n$, 
\[
\Phi^+_1(n) = d_1(n).
\]
In particular, $\Phi^+_1(n) = 1$ 
if and only if $n$ is a power of 2.
\et

Equivalently, the positive integer $n$ is a sum of $k \geq 2$ consecutive  positive integers 
if and only if $n$ is not a power of 2.  

\begin{proof}
Let $k$ be an odd  positive divisor of $n$.  
The identities
\[
a_{k,1}(n) = \frac{n}{k} - \frac{(k-1)}{2}
\qqand
b_{k,1}(n) = \frac{k + 1}{2} - \frac{n}{k}
\]
imply that 
\[
a_{k,1}(n)+b_{k,1}(n)=1
\]
and so exactly one of the integers $a_{k,1}(n)$ and $b_{k,1}(n)$ is positive.
Thus, for each odd positive divisor $k$ of $n$ there is exactly 
one sequence of consecutive positive integers that sums to $n$.
This proves that $\Phi^+_1(n) = d_1(n)$.
\end{proof}

\bt                     \label{trapezoid:theorem:t-odd}
For every  odd positive integer $t$,  
let 
\[
\theta_t(n) =  \sqrt{\frac{2n}{t}+\frac{1}{4}} + \frac{1}{2} 
\qqand
\psi_t(n) =  \sqrt{2nt + \left( \frac{t-2}{2} \right)^2  } -  \left( \frac{t-2}{2} \right). 
\]
The number of representations of $n$ as a positive trapezoid with difference $t$ is 
\beq   \label{trapezoid:Phi+}
\Phi_t^+(n) = d_1(n) + d_1(n, \theta_t(n)) - d_1(n, \psi_t(n)).
\eeq
\et

\begin{proof}
Let $k$ be an odd divisor of $n$.  
All of the summands in the length $k$ representation~\eqref{trapezoid:odd} 
are positive if and only if $a_{k,t}(n) > 0$, or, equivalently, $k < \theta_t(n)$.  
The number of such divisors is $d_1(n,\theta_t(n))$.  

All of the summands in the length $2n/k$ representation~\eqref{trapezoid:even} 
are positive if and only if $b_{k,t}(n) > 0$ or, equivalently, $k \geq \psi_t(n)$.
The number of such divisors is $d_1(n) - d_1(n,\psi_t(n))$.  
This completes the proof.
\end{proof}

Note that if $t=1$, then $\theta_1(n) = \psi_1(n)$, and so, 
for every odd divisor $k$ of $n$, exactly one of the inequalities 
$k < \theta_1(n)$ and $k \geq \psi_1(n)$ will hold.  
This gives another proof that $\Phi^+_1(n) = d_1(n)$.

In a \emph{Comptes Rendus} note in 1883, Sylvester~\cite{sylv83} proved that 
``\ldots le nombre de suites de nombres cons{\' e}cutifs dont la somme 
est \N\ est {\' e}gal au nombre de diviseurs impairs de \N.''
This result (Theorem~\ref{trapezoid:theorem:PowerOf2}) 
 has been rediscovered many times.  
A special case is in \emph{Number Theory for Beginners}~\cite{weil79b} 
by  Andr\' e Weil:  
Problem III.4  is to prove that an ``integer $> 1$ which is not a power of 2 can be written 
as the sum of 2 or more consecutive integers.''

MacMahon~\cite[vol. 2, p. 28]{macm60} used generating functions 
to prove Theorem~\ref{trapezoid:theorem:PowerOf2}.  

Here is a nice generalization.   
Let $\Phi^+_{1,0}(n)$ (resp. $\Phi^+_{1,1}(n)$) denote the number of representations of $n$ 
as the sum of an even (resp. odd) number of consecutive positive integers.  
Thus, 
\[
\Phi^+_1(n) = \Phi^+_{1,0}(n)  + \Phi^+_{1,1}(n). 
\]
Andrews, Jim\'enez-Urroz,  and Ono~\cite{andr01} proved analytically that 
\[
\Phi^+_{1,0}(n)  - \Phi^+_{1,1}(n) = d(n, \sqrt{2n}) - d(n, \sqrt{ (n/2)}).
\]
Chapman~\cite{chap02} gave a combinatorial proof of this result.

\section{Hook numbers and the Durfee square}
Before describing Sylvester's algorithm,  we recall some properties of the Durfee square 
of a partition of a positive integer $n$.  
Let 
\beq        \label{trapezoid:r-part}  
n = r_1 + \cdots + r_k
\eeq
be a partition of $n$ into $k$ positive and decreasing parts.   
We have $r_1 \geq 1$. 
Let $s$ be the greatest integer such that $r_s \geq s$.  
The square array of $s^2$ dots in the upper left corner of the Ferrers
graph is called the \emph{Durfee square}\index{Durfee square} of the partition, 
and the positive integer $s$ is the \emph{side} of the Durfee square.  
If $s+1 \leq i \leq k$, then $r_i \leq r_{s+1} \leq s$ and all of the dots on the $i$th row 
of the Ferrers graph lie on the first $s$ columns of the graph.  
It follows that every dot in the Ferrers graph lies on one of the first 
$s$ rows or on one of the first $s$ columns of the graph.  
Therefore, the row numbers $r_1,\ldots, r_s$ 
and the column numbers $c_1,\ldots, c_s$ determine the 
partition~\eqref{trapezoid:r-part}.  
We extend this observation as follows.

\bl                \label{trapezoid:lemma:rc-unique}
Let $s_1$ and $s_2$  be positive integers, and let $(r_i)_{i=1}^{s_1}$ and $(c_j)_{j=1}^{s_2}$ 
be sequences of integers such that 
\[
r_1 \geq r_2 \geq \cdots \geq r_{s_1} \geq s_2
\]
and 
\[
c_1 \geq c_2 \geq \cdots \geq c_{s_2} \geq s_1.
\]
The positive integer  
\[
n = \sum_{i=1}^{s_1} r_i  + \sum_{j =1}^{s_2} c_j  - s_1 s_2
\]
has a unique partition 
with parts $r_1,\ldots, r_{s_1}, r_{s_1 +1},\ldots, r_{c_1}$, 
where, for $i = s_1 +1, \ldots, c_1$, 
\[
r_i = \max(j : c_j \geq i).
\] 
If $s_1 = s_2 = s$, then the Durfee square of this partition has side $s$, 
and the row numbers $r_1,\ldots, r_s$ and column numbers 
$c_1,\ldots, c_s$ determine the partition.  
\el

\begin{proof}
Note that 
\[
n = \sum_{i=1}^{s_1} r_i  + \sum_{j =1}^{s_2} (c_j  - s_1) \geq  \sum_{i=1}^{s_1} r_i \geq s_1s_2.
\]
Construct the Ferrers diagram with $r_i$ dots on row $i$ for $i=1,\ldots, s_1$, 
and with $c_j$ dots on column $j$ for $j=1,\ldots, s_2$.  
The Ferrers diagram has $c_1$ rows, 
and so the partition of $n$ has $c_1$ parts.  
For $i = s_1 +1,\ldots, c_1$, 
there is a dot on the $j$th column of row $i$ if and only if $j \leq s_2$ and $c_j \geq i$.  
Therefore, $r_i =  \max(j : c_j \geq i) \leq s_2$.  

If $s_1 = s_2 = s$, then $r_{s+1} = r_{s_1+1} \leq s_2  \leq r_{s_1} = r_s$, 
and so this partition has a Durfee square with side $s$.  
This completes the proof.  
\end{proof}

The upper left corner of a Ferrers diagram of a partition 
contains a unique minimal square array of dots 
(the Durfee square) whose rows and columns determine the partition.  
The upper left corner of a Ferrers diagram also contains minimal rectangular arrays of dots 
whose rows and columns determine the partition. 
The Ferrers diagram contains a ``Durfee rectangle''\index{Durfee rectangle} 
with sides $(s_1,s_2)$ if 
\[
r_{s_1 + 1} \leq s_2 \leq r_{s_1} \qqand c_{s_2 + 1} \leq s_1 \leq c_{s_2}.
\]
These Durfee rectangles are not unique. For example, the partition 
\[
23 = 5 + 5 + 4 + 3 + 3 + 2 + 1
\]
has Durfee square of side 3, and Durfee rectangles of sides $(s_1,s_2) = (2,4)$ 
and $(s_1,s_2) = (5,2)$.

For $1 \leq i \leq k$ and $1 \leq j \leq r_i$, let $R_{i,j} $ be the set of dots 
on the $i$th row that are on and to the right of the $j$th dot, 
and let $C_{i,j} $ be the set of dots on the $j$th column 
that are on and below the $i$th dot.  
The \emph{$(i,j)$th hook number}\index{hook number} is the cardinality of the set 
$H_{i,j} = R_{i,j} \cup C_{i,j}$.
The number of dots on row $i$ is $r_i = |R_{i,1}|$.  
Denote the number of dots on 
column $j$ by $c_j = |C_{1,j}|$.  We obtain  
\[
|H_{i,j}| = r_i + c_j - i - j + 1.
\]

For $i=1,\ldots, s$, we define the 
\emph{diagonal hook number}\index{diagonal hook number} 
\[
h_i = |H_{i,i}| = r_i + c_i - 2i + 1.
\] 
The set of diagonal hooks $\{H_{i,i} : i = 1,\ldots, s\}$ partitions the dots 
in the Ferrers diagram 
and produces the \emph{hook partition}\index{hook partition} of $n$: 
\[
n = h_1 + h_2 + \cdots + h_s.  
\]

\bl              \label{trapezoid:lemma:h-diff}
Let $n = r_1 + \cdots + r_k$ be a partition of $n$,  let $s$ be the side of the Durfee square 
of the Ferrers diagram of this partition, and let $h = h_1 + \cdots + h_s$ 
be the associated hook partition of $n$.  For $i=1,\ldots, s-1$,
\[
h_i - h_{i+1} \geq 2
\]
and 
\[
h_i - h_{i+1} = 2
\]
if and only if $r_i = r_{i+1}$ and $c_i = c_{i+1}$.  
\el

\begin{proof}
For $i=1,\ldots, s-1$ we have 
\begin{align*}
h_i - h_{i+1} 
& = (  r_i + c_i - 2i + 1 ) -  (  r_{i+1} + c_{i+1} - 2i - 1 )    \\
& = (  r_i - r_{i+1}  ) + ( c_i - c_{i+1} ) + 2  \\
& \geq 2.
\end{align*}
Moreover, $h_i - h_{i+1} = 2$ if and only if $r_i = r_{i+1}$ and $c_i = c_{i+1}$.  
\end{proof}

For example, the partition into odd parts
\[
57 = 11 + 11 + 11 + 9 + 5 + 5 + 5
\]
has the left-justified Ferrers graph 
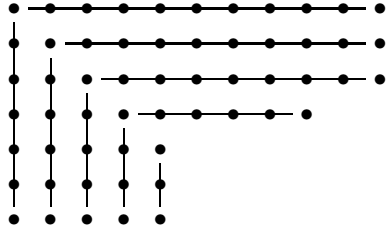
\[
\xymatrix@=0.1cm{
\bullet \ar@{-}[dddddd] \ar@{-}[rrrrrrrrrr] & \bullet & \bullet  & \bullet & \bullet & \bullet & \bullet & \bullet &\bullet & \bullet & \bullet  \\
\bullet  & \bullet \ar@{-}[ddddd]  \ar@{-}[rrrrrrrrr]  & \bullet & \bullet & \bullet & \bullet & \bullet & \bullet &\bullet & \bullet & \bullet \\
\bullet  & \bullet  & \bullet \ar@{-}[dddd]  \ar@{-}[rrrrrrrr]  & \bullet & \bullet & \bullet & \bullet & \bullet &\bullet & \bullet & \bullet \\
\bullet  & \bullet  & \bullet & \bullet \ar@{-}[ddd] \ar@{-}[rrrrr] & \bullet & \bullet & \bullet  & \bullet &\bullet  \\
\bullet  & \bullet  & \bullet & \bullet & \bullet \ar@{-}[dd]\\
\bullet  & \bullet  & \bullet & \bullet & \bullet \\
\bullet  & \bullet & \bullet & \bullet & \bullet
}
\]
We have $5 = r_5 = r_6 < 6$ and so  the Durfee square has side 5 
contains $5^2 = 25$ dots.  The hook partition is
\[
57 = 17 + 15 + 13 + 9 + 3.
\]

Note that the hook partition of a partition does not determine the partition.
For example, the partitions $5 + 2$ and $4 + 2 + 1$ both have Durfee squares 
of side 2 and hook partitions $6 + 1$.  
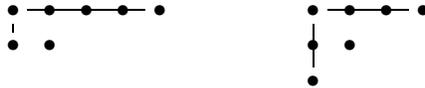
\[
\xymatrix@=0.1cm{
\bullet \ar@{-}[d] \ar@{-}[rrrr] & \bullet & \bullet  & \bullet & \bullet &&&&&&  \bullet \ar@{-}[dd] \ar@{-}[rrr] & \bullet & \bullet  & \bullet  \\
\bullet  & \bullet &&&&&&&&& \bullet & \bullet \\
&&&&&&&&&& \bullet &  \\
}
\]

\bt
The number of partitions of $n$ into exactly $k$ parts differing by at least 2 is the 
number of partitions of $n-k^2$ into at most $k$ parts.  
\et

\begin{proof}
The first construction converts a partition of $n-k^2$ into at most $k$ parts
into a partition of $n$ into exactly $k$ parts differing by at least 2. 
Let $n > k^2$, and let $n-k^2 = \sum_{i=1}^k b_i$ be a partition 
with $1 \leq r \leq k$ and $b_1 \geq \cdots \geq b_r$. 
For $r+1 \leq i \leq k$ we define $b_i = 0$, and for $i=1,\ldots, k$ we define 
\[
a_i = b_i + 2(k-i) +1.
\] 
It follows that 
\begin{align*}
a_i  - a_{i+1} & =  (b_i + 2(k-i) +1) - (b_{i+1} + 2(k-i-1) +1) \\ 
& =  b_i - b_{i+1} + 2 \geq 0
\end{align*}
for $i=1,\ldots, k-1$.  The identity 
\[
k^2 = \sum_{i=1}^k (2i-1) =  \sum_{i=1}^k (2(k-i) +1) 
\]
implies that 
\[
n = \left(n-k^2\right) + k^2 = \sum_{i=1}^k (b_i + 2(k-i) +1 ) = \sum_{i=1}^k a_i.
\]
This is a partition of $n$ into exactly $k$ parts differing by at least 2.  

The second construction converts a partition of $n$ into exactly $k$ parts 
differing by at least 2 into a partition of $n-k^2$ into at most $k$ parts.  
Let $n = \sum_{i=1}^k a_i$ be a partition of $n$ into exactly $k$ parts 
differing by at least 2.  We have $a_k \geq 1$.  
If $1 \leq i \leq k-1$ and $a_{i+1} \geq 2(k-(i+1)) + 1$, 
then 
\[
a_i \geq a_{i+1}+2 \geq (2(k-(i+1)) + 1) + 2 = 2(k- i) + 1.
\]
It follows by downward induction that $a_i \geq 2(k- i) + 1$ and so   
\[
b_i = a_i - (2(k- i) + 1) \geq 0
\]
for $i=1,\ldots, k$.   
We have  
\[
\sum_{i=1}^k b_i  = \sum_{i=1}^k a_i - \sum_{i=1}^k (2(k- i) + 1) = n - k^2.
\]
This is a partition of $n-k^2$ into at most $k$ parts.  

It is straightforward to check that the first and second constructions are inverses of each other.  
This completes the proof.   
\end{proof}

Consider a partition of $n$ whose Ferrers diagram has Durfee square of side $s$.    
Let $r_1,\ldots, r_s$ be the number of dots on the first $s$ rows of the Ferrers diagram, 
and let $c_1,\ldots, c_s$ be the number of dots on the first $s$ columns.  
The \emph{Frobenius symbol}\index{Frobenius symbol} of the partition is the $2 \times s$ matrix
\[
\bmat r_1 - 1 & r_2 - 2 & \cdots & r_s - s \\
 c_1 - 1 & c_2 - 2 & \cdots & c_s - s 
 \emat.
\]
Note the rows are strictly decreasing sequences of nonnegative integers, 
and that 
\[
n = s + \sum_{i=1}^s (r_i-1) + \sum_{i=1}^s (c_i-1).  
\]
The Frobenius symbol is related to the construction in Lemma~\ref{trapezoid:lemma:rc-unique}.  
See Andrews~\cite{andr84a,andr84b}.

\section{Sylvester's algorithm}    \label{trapezoid:section:Sylvester}
Sylvester discovered a graphical algorithm, sometimes called the 
\emph{fish-hook method}, that transforms a partition of $n$ with odd parts 
into a partition of $n$ with distinct parts, and showed that this transformation is a bijection 
between the set of partitions into odd parts and the set of partitions into distinct parts.  
Moreover, he proved that this transformation has the extraordinary  property that 
if the original partition of $n$ into odd parts 
contains exactly $\ell$ different odd integers, then the new partition of $n$ into 
distinct parts contains exactly $\ell$ maximal subsequences of consecutive integers.

Here is the algorithm.  
Let 
\beq                                  \label{trapezoid:SA1}
n = a_1 + \cdots + a_k
\eeq
be a partition of $n$ into odd parts, with 
\beq                                  \label{trapezoid:SA2}
a_1 \geq \cdots \geq a_k \geq 1
\eeq 
and  
\beq                                  \label{trapezoid:SA3}
a_i = 2r_i - 1
\eeq
for $i=1,\ldots, k$.  Then
\beq                                  \label{trapezoid:SA4}
r_1 \geq \cdots \geq r_k \geq 1.
\eeq
Because the summands $a_i$ are odd, we can draw a center-justified Ferrers diagram,
and divide it into two sub-diagrams.
The \emph{major right half} consists of the vertical central line and the dots to its right.
The \emph{minor left half} consists of the dots that are strictly to the left of the central line.
We compute the hook numbers of the major half, and denote them in decreasing order by 
$h_1 > h_3 > h_5 > \cdots.$ 
We compute the hook numbers of the minor half, and denote them in decreasing order by 
$h_2 > h_4 > h_6 > \cdots.$ 
We shall prove that $h_1 > h_2 > h_3 > h_4 > h_5 > h_6 \cdots$, 
and so the hook numbers create a partition of $n$ into distinct parts.

Before proving this statement, we consider an example: 
\[
57 = 11 + 11 + 11 + 9 + 5 + 5 + 5
\]
is a partition into odd parts.  
The center-justified Ferrers diagram is 
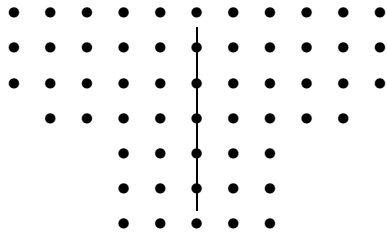
\[
\xymatrix@=0.1cm{
\bullet & \bullet & \bullet & \bullet & \bullet & \bullet \ar@{-}[dddddd] & \bullet & \bullet &\bullet & \bullet & \bullet \\
\bullet & \bullet & \bullet & \bullet & \bullet & \bullet & \bullet & \bullet &\bullet & \bullet & \bullet \\
\bullet & \bullet & \bullet & \bullet & \bullet &\bullet & \bullet & \bullet &\bullet & \bullet & \bullet \\
& \bullet & \bullet & \bullet & \bullet & \bullet & \bullet & \bullet & \bullet &\bullet & \\
& & & \bullet & \bullet &\bullet & \bullet & \bullet & & &\\
& & & \bullet & \bullet & \bullet& \bullet & \bullet & & &\\
& & & \bullet & \bullet & \bullet & \bullet & \bullet & & &
}
\]  
The major half is the  Ferrers diagram 
of the partition $32 = 6 + 6 + 6 + 5 + 3 + 3 + 3$:
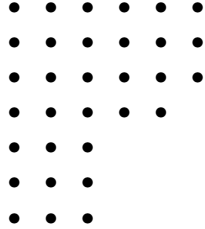
\[
\xymatrix@=0.1cm{
 \bullet & \bullet & \bullet &\bullet & \bullet & \bullet \\
\bullet & \bullet & \bullet &\bullet & \bullet & \bullet \\
\bullet & \bullet & \bullet &\bullet & \bullet & \bullet \\
\bullet & \bullet & \bullet & \bullet &\bullet & \\
\bullet & \bullet & \bullet & & &\\
 \bullet& \bullet & \bullet & & &\\
 \bullet & \bullet & \bullet & & &
}
\]
The remainder of the original Ferrers diagram is the minor half, 
associated with the partition $25 = 5 + 5 + 5 + 4 + 2 + 2 + 2$:
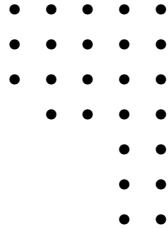
\[
\xymatrix@=0.1cm{
\bullet & \bullet & \bullet & \bullet & \bullet  \\
\bullet & \bullet & \bullet & \bullet & \bullet \\
\bullet & \bullet & \bullet & \bullet & \bullet \\
          & \bullet & \bullet & \bullet & \bullet  \\
                     & & & \bullet &\bullet \\
                    & & & \bullet & \bullet \\
                   & & & \bullet & \bullet
}
\]
which we rearrange as the Ferrers diagram 
of the partition $5 + 5 + 5 + 4 + 2 + 2 + 2$:
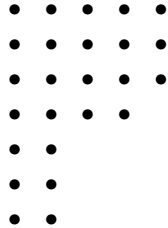
\[
\xymatrix@=0.1cm{
\bullet & \bullet & \bullet & \bullet & \bullet  \\
\bullet & \bullet & \bullet & \bullet & \bullet \\
\bullet & \bullet & \bullet & \bullet & \bullet \\
         \bullet & \bullet & \bullet & \bullet   & \\
                      \bullet &\bullet & & &\\
                    \bullet & \bullet  & & &\\
                    \bullet & \bullet & & &
}
\]
Note that deleting the first column of the major half produces the minor half.

The Durfee square of the major half consists of $4^2 = 16$ vertices.  
Every dot in this diagram lies on one of the first four rows 
or on  one of the first four columns.  
We partition the vertices  of the major  half into the four hooks  of the Durfee square  
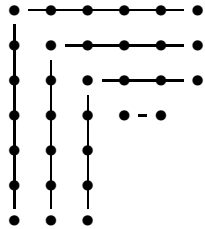
\[
\xymatrix@=0.1cm{
 \bullet \ar@{-}[rrrrr] \ar@{-}[dddddd]  & \bullet & \bullet &\bullet & \bullet & \bullet  \\
\bullet  & \bullet  \ar@{-}[rrrr] \ar@{-}[ddddd]  & \bullet &\bullet & \bullet & \bullet \\
\bullet  & \bullet  & \bullet   \ar@{-}[rrr] \ar@{-}[dddd] &\bullet & \bullet & \bullet \\
\bullet  & \bullet  & \bullet  & \bullet  \ar@{-}[r]  &\bullet & \\
\bullet  & \bullet  & \bullet  & & &\\
 \bullet & \bullet  & \bullet  & & &\\
 \bullet & \bullet & \bullet & & &
}
\]
and obtain the hook partition   
\[
32 = 12 + 10 + 8 + 2.
\]

The minor left half is the major half with the left column removed,  
and the Durfee square of the minor half also consists of $16$ vertices. 
Separating the minor half into hooks, we obtain 
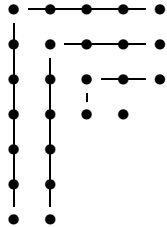
\[
\xymatrix@=0.1cm{
\bullet \ar@{-}[dddddd]  \ar@{-}[rrrr] & \bullet  & \bullet  & \bullet  & \bullet   \\
\bullet  & \bullet \ar@{-}[ddddd] \ar@{-}[rrr]  & \bullet  & \bullet  & \bullet  \\
\bullet  & \bullet  & \bullet  \ar@{-}[d] \ar@{-}[rr] & \bullet  & \bullet  \\
           \bullet  & \bullet & \bullet  & \bullet   & \\
                     \bullet  &\bullet  & & & \\
                     \bullet  & \bullet  & & & \\
                   \bullet & \bullet& & & 
}
\]
with hook partition 
\[
25 = 11 + 9 + 4 + 1.
\]
Notice that not only are the parts in the hook partitions strictly decreasing, 
but they are also interlaced in magnitude.  
Their union gives a partition of 57 into distinct parts:
\[
57 = 12 + 11 +  10 + 9 + 8 + 4 + 2 + 1.
\]
Thus, the original partition with odd parts has been transformed 
into a partition with distinct parts.  
We also observe that the original partition of 57 used 
only the three odd integers 11, 9, and 5,
and that the new partition of 57 into distinct parts consists 
of three maximal decreasing sequences of 
consecutive integers:  $(12,11,10,9,8)$, $(4)$, and $(2,1)$.  

MacMahon~\cite[vol. 2, pp. 13--14]{macm60} 
contains a description of Sylvester's fish-hook method.  
Andrews~\cite[Section 4]{andr84a} uses the Frobenius symbol of a partition 
to explain the fish-hook method.

\section{Sylvester's proof of Euler's theorem}
\bt              \label{trapezoid:theorem:Sylvester}
Let $n$ be a positive integer, 
let $\mcu(n)$ be  the set of all partitions of $n$ into odd parts, 
and $\mcv(n)$ be the set of all partitions of $n$ into distinct parts.  
The function $f:\mcu(n) \rightarrow \mcv(n)$ defined by Sylvester's algorithm is a bijection.  
\et

\begin{proof}
Consider a partition of $n$ into $k$ odd parts of the 
form~\eqref{trapezoid:SA1} -- ~\eqref{trapezoid:SA4}.
Let $s$ be the side of the Durfee square of the major half.   
Every dot in the major half  lies on one of the first $s$ rows 
or on one of the first $s$ columns.
For $i=1,\ldots, s$, the number $r_i$ of dots on the $i$th row of the major half 
satisfies
\[
r_1 \geq \cdots \geq r_s \geq s \geq r_{s+1}.  
\]
%Moreover, $r_i = 1$ if and only if $a_i = 1$.

For $i=1,\ldots, r_1$,  let 
$c_i$ be the number of dots in the $i$th column of the major half.
Note that $r_{s+1} \leq s$ implies that $c_{s+1} \leq s$, and so   
\[
k = c_1 \geq \cdots \geq c_s \geq s \geq c_{s+1}.
\]
For $i=1,\ldots, s$, we have the hook numbers
\beq   \label{trapezoid:RightHook}
h_i = r_i + c_i - 2i + 1.
\eeq
By Lemma~\ref{trapezoid:lemma:h-diff}, these numbers 
satisfy $h_i - h_{i+1} \geq 2$ for $i=1,\ldots, s -1$.

The minor half of the original Ferrers diagram is exactly the major half 
with the first column removed.  
Therefore, every dot in the minor half lies on one of the first $s$ rows of the  
minor half or on one of the first $s-1$ columns of the graph.  
For $i = 1,\ldots, c_2$, let $r'_i= r_i - 1$ denote the number of dots 
on the $i$th row of the minor half.  
For $i = 1,\ldots, r'_1$, let  $c'_i$ denote the number of dots 
on the $i$th column of the minor half.

%Moreover,  $r'_i = 0$ if and only if $a_i = 1$,  and $c'_s = 0$ if and only if $r_1 = s$.  

Let $s'$ be the side of the Durfee square of the minor half.  
Because
\[
r'_{s-1} \geq r'_s = r_s - 1 \geq s-1 \geq r_{s+1} -1 = r'_{s+1}
\]
it follows that $s' = s-1$ or $s' = s$,  
Moreover,  $s'=s$ if and only if $r_s \geq s + 1$ and $c_{s+1} = s$.  
Similarly, $s' = s-1$ if and only if $r_s = s$ and $c_{s+1} = s-1$.  .   

For $i=1,\ldots, s'$, there are the hook numbers
\beq   \label{trapezoid:LeftHook}
h'_i = r'_i + c'_i - 2i + 1 = r_i + c_{i+1}-2i.
\eeq
By~Lemma~\ref{trapezoid:lemma:h-diff}, we have $h'_i - h'_{i+1} \geq 2$ for $i=1,\ldots, s'-1$.

If $s' = s$, then $r_s \geq s+1$ and $c_{s+1} = s$,
and so
\[
h'_s = r_s + c_{s+1}-2s = r_s - s.
\] 
If $s' = s-1$, then $r_s  = s$.  
We define $h'_s = 0$, and again have
\[
h'_s = r_s - s
\]
and 
\begin{align*}
h'_{s-1} - h'_s & = h'_{s-1} =  r_{s-1} + c_s -2s+ 2 \\
& \geq (r_s-s) + (c_s - s) + 2 \\ 
& \geq 2.
\end{align*}

%Note that $c'_s = c_{s+1} = s$ if and only if $r'_s = r_s -1\geq s$.

We shall prove that 
\[
h_1 > h'_1 > h_2 > h'_2 > \cdots > h_s > h'_s \geq 0.
\]
For $i=1,\ldots, s-1$, we have
\begin{align*}
h_i - h'_i & = (r_i + c_i - 2i + 1) - (r'_i + c'_i - 2i + 1) \\ 
& = (r_i-r'_i)+(c_i-c'_i) \\ 
& = 1 + c_i-c_{i+1} \\ 
& \geq 1.
\end{align*}
Also,
\begin{align*}
h_s - h'_s & = (r_s+c_s-2s+1) - (r_s-s) \\
& = c_s-s+1 \geq 1.
\end{align*}
For $i=1,\ldots, s-1$, we have  
\begin{align*}
h'_i - h_{i+1} & = (r'_i + c'_i - 2i + 1) -  (r_{i+1} + c_{i+1} - 2i - 1) \\
& = (r_i -1 + c_{i+1} - 2i + 1) -  (r_{i+1} + c_{i+1} - 2i - 1) \\
& = r_i - r_{i+1} +1  \\
& \geq 1.
\end{align*}
Therefore, 
\beq         \label{trapezoid:distinct-even}
n = h_1 +  h'_1 + h_2 + h'_2 + \cdots + h_{s-1} +  h'_{s-1} + h_s + h'_{s}
\eeq
is a partition into $2s$ or $2s-1$ distinct positive parts, 
and we have transformed a  partition with only odd parts to a 
partition into distinct parts.
We shall prove that this transformation is one-to-one and onto.  

Consider a partition of $n $ into $2s$ distinct nonnegative parts:  
\[
n = h_1 + h'_1 + h_2 + h'_2 + \cdots + h_{s} + h'_{s}
\]
where 
\[
h_1 >  h'_1 > h_2 > h'_2 > \cdots > h_{s-1} >  h'_{s-1} > h_s > h'_s  \geq 0.
\]
If the number of positive parts is even, then $h'_s \geq 1$.  
If the number of positive parts is odd, then $h'_s = 0$.  

If this partition is constructed  by Sylvester's algorithm  
from a partition of $n$ into odd parts,  
then there are positive integers $r_1, r_2, \ldots, r_s$ and $c_1, c_2,\ldots, c_s$ 
such that 
\begin{align*}
h_1 & = r_1 +c_1 - 1   \\
h'_1 & = r_1 +c_2 - 2   \\
& \vdots  \\
h_i & = r_i  +c_i  - (2i -1)  \\
h'_i & = r_i  +c_{i + 1}  - 2i     \\
& \vdots  \\
%h_{s-1}  & = r_{s-1}  +c_{s-1}  - (2s-3)   \\
%h'_{s-1} & = r_{s-1} +c_{s-1} - (2s-2)   \\
h_s & = r_s +c_s - (2s-1)   \\
h'_s & = r_s - s.
\end{align*}
Conversely, given the $2s$ parts $h_1, h'_1, \ldots, h'_s$, 
we can solve these $2s$ equations recursively, and obtain unique 
integers $r_1,\ldots, r_s, c_1,\ldots, c_s$.  
For $i=1,\ldots, s$, the inequality $h_i > h'_i $ implies that 
\[
r_i  +c_i  - (2i -1) >  r_i  +c_{i + 1}  - 2i  
\]
and so
\[
c_i \geq c_{i+1}
\]
For $i=1,\ldots, s-1$, the inequality $h'_i > h_{i+1} $ implies that 
\[
 r_i  + c_{i + 1}  - 2i  >  r_{i+1}  + c_{i+1}  - (2i +1) 
\]
and so
\[
r_i \geq r_{i+1}.  
\]
Because 
\[
r_s = h'_s + s \geq s
\]
and
\begin{align*}
c_s & = h_s - r_s + 2s-1 \\
& = h_s - (h'_s + s) + 2s-1 \\
& = h_s - h'_s + s-1 \\
& \geq s
\end{align*}
it follows that $r_1 \geq \cdots \geq r_s \geq s$ and $c_1 \geq \cdots \geq c_s \geq s$ 
are decreasing sequences of positive integers.  

Thus, every partition into odd parts determines a unique partition into distinct parts, 
and every partition into distinct parts can be obtained uniquely from a partition 
into odd parts.  
\end{proof}

For example, consider the partition 
\[
50 = 22 + 17 + 8 + 3.
\]
We have $d=2$ and 
\begin{align*}
22 & = r_1 +c_1 - 1   \\
17 & = r_1 +c_2 - 2   \\
8 & = r_2 +c_2 - 3   \\
3 & = r_2-2.
\end{align*}
Solving these equations, we obtain 
\begin{align*}
r_2 & = 5  \\
c_2 & = 6  \\
r_1 & = 13  \\
c_1 & = 10.
\end{align*}
Thus, the major half has 10 rows, of lengths
\begin{align*}
r_1 & = 13  \\
r_2 & = 5  \\
r_i & = 2 \quad \text{for $i =  3, \ldots, 6$ }\\
r_i & = 1 \quad \text{for $i =  7, \ldots, 10$.}
\end{align*}
Defining $a_i = 2r_i-1$ for $i=1,\ldots, 10$, we obtain 
the following partition of 50 into odd parts:
\[
50 = 25 + 9 + 3 + 3 + 3 + 3 + 1 + 1 + 1 + 1.
\]
Note that the partition into distinct parts consists of four maximal sequences 
of consecutive integers, and that the corresponding partition into odd parts 
contains four distinct odd numbers.  

Here is another example:
\[
31 = 9 + 8 + 7 + 4 + 3.
\]
We have $d=3$ and 
\begin{align*}
9 & = r_1 +c_1 - 1   \\
8 & = r_1 +c_2 - 2   \\
7 & = r_2 +c_2 - 3   \\
4 & = r_2 + c_3 - 4\\
3 & = r_3  + c_3 - 5 \\
0 & = r_3 - 3.
\end{align*}
Solving these equations, we obtain 
\begin{align*}
r_3 & = 3 \\
c_3 & = 5 \\
r_2 & = 3  \\
c_2 & = 7  \\
r_1 & = 3  \\
c_1 & = 7.
\end{align*}
Thus, the major half has 7 rows, of lengths
\begin{align*}
r_i & = 3 \quad \text{for $i =  1, 2, 3, 4$ }\\
r_i & = 2 \quad \text{for $i =  5,6,7$.}
\end{align*}
Defining $a_i = 2r_i-1$ for $i=1,\ldots, 7$, we obtain 
the following partition of 31 into odd parts:
\[
31 =  5 + 5 + 5 + 5 + 5 + 3 + 3.
\]
Note that the partition into distinct parts consists of two maximal sequences 
of consecutive integers, and that the corresponding partition into odd parts 
contains two distinct odd numbers.  

Another example:
The partition into distinct parts 
\[
30 = 10 + 8 + 7 + 4 + 1 
\]
is mapped to the following partition into odd parts:
\[
30 = 9 + 9 + 5 + 3 + 3 + 1.
\]

\section{Sylvester's stratification of Euler's theorem}

For every positive integer $n$, let $p_{\text{odd}}(n) = |\mcu(n)|$, 
where $\mcu(n)$ is the set of partitions of $n$ into not necessarily distinct odd parts.
Let $p_{\text{dis}}(n) = |\mcv(n)|$, where $\mcv(n)$ is the set of partitions 
of $n$ into distinct parts.  
Euler proved (Theorem~\ref{trapezoid:theorem:Euler-odd-distinct}) 
that these two sets have the same cardinality, that is, $p_{\text{odd}}(n) = p_{\text{dis}}(n)$.  
In the proof of Theorem~\ref{trapezoid:theorem:Sylvester}, 
we proved that the function $f:\mcu(n) \rightarrow \mcv(n)$ 
defined by Sylvester's algorithm is a bijection.

For positive integers $n$ and $\ell$, let $\mcu_{\ell}(n)$ denote the set 
of partitions of $n$ into not necessarily distinct odd parts 
with exactly $\ell$ distinct odd parts, and let $U_{\ell}(n) = |\mcu_{\ell}(n)|$.  
We have 
\[
p_{\text{odd}}(n) = \sum_{\ell = 1}^{\infty} U_{\ell}(n).
\]
Similarly, if $\mcv_{\ell}(n)$ denotes the set 
of partitions of $n$ into distinct parts and $V_{\ell}(n) = |\mcv_{\ell}(n)|$, then 
\[
p_{\text{dis}}(n) = \sum_{\ell = 1}^{\infty} V_{\ell}(n).
\]
Sylvester's ``stratification'' of Euler's theorem is that $U_{\ell}(n) = V_{\ell}(n)$ 
for all positive integers $n$ and $\ell$.  

For example, the set $\mcu_3(57)$ contains the partition 
\[
11 + 11 + 11 + 9 + 5 + 5 + 5
\]
which is a partition of 57 into odd parts whose three distinct parts are 11, 9, and 5.  
Similarly, the set $\mcv_3(57)$ contains the partition 
\[
12 + 11 + 10 + 9 + 8 + 4 + 2 + 1
\]
which is a partition of 57 with three maximal subsequences 
of consecutive integers:  \\
$(12, 11, 10, 9, 8)$, $(4)$, and $(2,1)$.

There are three partitions of 5 into odd parts: $5 = 3+1+1 = 1+1+1+1+1$.  
The partitions with one distinct part are $5$ and $1+1+1+1+1$, 
and so $U_1(5) = 2$.  The partition with two distinct parts is $3+1+1$, 
and so $U_2(5) = 1$.

There are three partitions of 5 into distinct parts: $5 = 4+1 = 3+2$.
The partitions with one maximal subsequence  
of consecutive integers are $5$ and $3+2$, and so $V_1(5) = 2$.   
The partition with two maximal subsequences  
of consecutive integers is $4+1$, and so $V_2(5) = 1$. 

The proof of Sylvester's theorem uses the following combinatorial observation.

\bl                  \label{trapezoid:lemma:counting}
Let \mcu\ and \mcv\ be sets, and let $\{\mcu_i: i=1,2,3,\ldots\}$ and 
$\{\mcv_i: i=1,2,3,\ldots\}$ be partitions of \mcu\ and \mcv, respectively.
Let $f:\mcu \rightarrow \mcv$ be a bijection.  
For every positive integer $\ell$, let $f_{\ell} :\mcu_{\ell} \rightarrow \mcv$ 
be the restriction of $f$ to $\mcu_{\ell}$. 
If $f_{\ell}(\mcu_{\ell}) \subseteq \mcv_{\ell}$ for all ${\ell} \in \N$, then 
$f_{\ell}:\mcu_{\ell} \rightarrow \mcv_{\ell}$ is a bijection for all ${\ell} \in \N$.
\el

\begin{proof}
Because $f$ is a bijection, it follows that $f$ is one-to-one, 
and so $f_{\ell}$ is one-to-one for all ${\ell} \in \N$.
Let $v\in \mcv_{\ell} \subseteq \mcv$.  
Because $f$ is onto, there exists $u\in \mcu$ such that $f(u) = v$.
Because $\mcu = \bigcup_{i=1}^{\infty} \mcu_i$ is a partition of \mcu, there is a unique 
integer $j$ such that $u \in \mcu_j$.  Therefore, $v = f(u) = f_j(u) \in \mcv_j$
and so $v \in \mcv_{\ell} \cap \mcv_j$.  
Because $\mcv = \bigcup_{i=1}^{\infty}  \mcv_i$ is a partition of \mcv, it follows that 
${\ell}=j$ and $u \in \mcu_{\ell}$.  Therefore, $f_{\ell}:\mcu_{\ell} \rightarrow \mcv_{\ell}$ 
is one-to-one and onto.
This completes the proof.  
\end{proof}

\bt                 \label{trapezoid:theorem:Sylvester-big}
Let 
\beq    \label{trapezoid:UsualPartition}
n = a_1 + \cdots + a_k
\eeq
be a partition of $n$ into $k$ not necessarily distinct odd  parts, 
and let $\ell$ be the number of distinct odd parts in this partition. 
The major-minor hook partition consists of exactly $\ell$ pairwise 
disjoint maximal sequences of consecutive integers.  
\et

\begin{proof}
Let 
\[
a_1 \geq a_2 \geq \cdots \geq a_k \geq 1
\]
and, for $i=1,\ldots, k$, let 
\[
a_i = 2r_i -1.
\]
We have   
\[
r_1 \geq r_2 \geq \cdots \geq r_k \geq 1.
\]
Let $\ell$ be the number of distinct odd parts in the partition~\eqref{trapezoid:UsualPartition}.  
The proof is by induction on $\ell$.  

If $\ell=1$, then $a_i = a_1 = 2r_1-1$ for $i=1,\ldots, k$, and $n = ka_1$.  
The Ferrers diagram for the partition is a rectangular array 
consisting of $k$ rows of $a_1$ dots.  
The major half of the diagram is a rectangular array 
consisting of $k$ rows of $r_1$ dots, and the minor half is 
a rectangular array consisting of $k$ rows of $r_1-1$ dots.  
The Durfee square of the major half has side $s = \min(k,r_1)$ and 
the Durfee square of the minor half has side $s' = \min(k,r_1-1)$.  
Let $n = h_1 + h'_1 + h_2 + h'_2 + \cdots$ be the major-minor hook partition.
By Lemma~\ref{trapezoid:lemma:h-diff}, we have 
\[
h_i - h_{i+1} = 2
\]
for $i=1,\ldots, s-1$, and
\[
h'_i - h'_{i+1} = 2
\]
for $i=1,\ldots, s'-1$.
Because 
\[
h_1 - h'_1 = (r_1+c_1-1) - (r_1+c_1-2) = 1 
\]
it follows that the parts in the major-minor hook partition of $n$ 
form a strictly decreasing sequence of consecutive integers.

For example, if $n = 21 = 7+7+7$, then $k = 3$, $r_1=4$,  
and the major-minor hook partition is $21 = 6+5+4+3+2+1$.  
If $n = 21 = 3 + 3 + 3 + 3 + 3 + 3 + 3$, then $k = 7$, $r_1=2$, 
and the major-minor hook partition is $21 = 8 +7 + 6$.

Let $\ell \geq 2$, and assume that the Theorem is true for partitions 
into at most $\ell -1$ distinct odd parts.  
The smallest part in the partition~\eqref{trapezoid:UsualPartition} 
is $a_k = 2r_k - 1$.  We also know that $ a_k < a_1 $ because $\ell \geq 2$.
If $j$ is the greatest integer such that $a_k < a_j $, 
then 
\begin{align*}
1 & \leq r_k   < r_j   \\
a_i & = a_k  \quad\text{ for } i=j+1,\ldots, k 
\end{align*}
and
\beq                \label{trapezoid:mg} 
m = n - (k-j)a_k = a_1 + \cdots + a_j
\eeq
is a partition of $m$ into odd parts with exactly $\ell - 1$ distinct parts.
By the induction hypothesis,  the Theorem is true for this partition of $m$.  

There are three cases.

Case 1:  
\[
\boxed{j < r_k < r_j}
\]
Because $j < r_j$, both the major  and the minor halves of the partition of $m$ have  
Durfee squares with side $j$. 
Let 
\beq                \label{trapezoid:Case1mg} 
m  = g_1 + g'_1 + \cdots + g_j + g'_j
\eeq
be the major-minor hook partition for $m$, where 
\beq               \label{trapezoid:Case1c-g}  
g_1 > g'_1 >g_2 > \cdots > g_j > g'_j.  
\eeq
For $i=1,\ldots, j$ we have
\begin{align}
g_i & = r_i + j - 2i + 1               \label{trapezoid:Case1a-g}    \\
g'_i & = r_i + j - 2i.      \label{trapezoid:Case1b-g} 
\end{align}
The partition~\eqref{trapezoid:mg} is  a partition of $m$ into odd parts with exactly $\ell - 1$ 
distinct parts.  By the induction hypothesis, the major-minor hook 
partition~\eqref{trapezoid:Case1mg} 
consists of exactly $\ell - 1$ pairwise disjoint maximal sequences of consecutive integers.  

Because
\[
j +1 \leq r_k = r_{j+1}
\]
the Durfee square for the major half of the partition of $n$ has side $s = \min(k,r_k) \geq j+1$, 
and the Durfee square for the minor half of the partition of $n$ 
has side $s' = \min(k,r_k-1) \geq j$.
Let 
\beq                \label{trapezoid:Case1nh} 
n  = h_1 + h'_1 + \cdots + h_j + h'_j + h_{j+1} + \cdots 
\eeq
be the major-minor hook partition for  $n$, 
where 
\[
h_1 > h'_1 >h_2 > \cdots > h_j > h'_j > h_{j+1} >\cdots. 
\]
For $i=1,\ldots, j$ we have
\begin{align}
h_i & = r_i + k - 2i + 1 = g_i + (k-j)        \label{trapezoid:Case1a}  \\
h'_i & = r_i + k - 2i = g'_i + (k-j)                 \label{trapezoid:Case1b}  
\end{align}
It follows that the number of pairwise disjoint maximal sequences of consecutive 
integers in the sequence $(g_1,\ldots, g'_1,\ldots, g_j,g'_j)$ 
of parts in the major-minor hook partition for $m$  is equal to the 
number of pairwise disjoint maximal sequences of consecutive integers 
in the sequence $(h_1,\ldots, h'_1,\ldots, h_j,h'_j)$.  
For $i= j+1,\ldots, s$ we have
\[
h_i  = r_k+ k - 2i  + 1 
\]
and for $i= j+1,\ldots, s'$ we have 
\[
h'_i  = r_k+ k - 2i. 
\]
We observe that, for $i > j$,  
\[
h_i - h'_i = h'_i - h_{i+1} = 1
\]
and so 
\beq         \label{trapezoid:case1}
(h_{j+1}, h'_{j+1}, h_{j+2}, \ldots)
\eeq
is a sequence of consecutive integers.  
Moreover, 
\[
h'_j - h_{j+1} = (r_j + k - 2j) - (r_k+ k - 2j -1) = r_j - r_k + 1 \geq 2
\]
and so~\eqref{trapezoid:case1} is a maximal  sequence of consecutive integers 
in the major-minor hook partition of $n$.  
It follows that the number of pairwise disjoint maximal sequences of consecutive 
integers in the major-minor hook partition of $n$ is exactly one more than the 
number of pairwise disjoint maximal sequences of consecutive 
integers in the major-minor hook partition of $m$.  
By the induction hypothesis, the latter partition consists of 
$\ell - 1$ maximal disjoint sequences, and so the partition~\eqref{trapezoid:Case1nh}
consists of $\ell$ maximal disjoint sequences.

For example, if 
\[
n = 49 = 13 + 13 + 9 + 7 + 7
\]
then $k=5$, $\ell = 3$, and 
\[
j = 3 < r_k = 4 < r_j = 5.
\]
We have $k-j = 2$ and 
\[
m = 35 = 13 + 13 + 9.  
\]
The major-minor hook partition for $m$ is 
\[
m = 35 = 9 + 8 + 7 + 6 + 3 + 2 
\]
and contains two maximal sequences of consecutive integers:
$(9,8,7,6)$ and $(3,2)$. 
The major-minor hook partition for $n$ is 
\[
n = 49 = 11 + 10 + 9 + 8 + 5 + 4 + 2
\]
and contains three maximal sequences of consecutive integers:
$(11, 10, 9,8)$, $(5,4)$, and $(2)$. 
Note that $(11, 10, 9,8) = (9,8,7,6) + (2,2,2,2)$ and $(5,4) = (3,2) + (2,2)$.  

Case 2: 
\[
\boxed{r_k \leq j < r_j}
\] 
Because $j < r_j$,  the major-minor hook partition of $m$ 
satisfies the relations~\eqref{trapezoid:Case1mg},
~\eqref{trapezoid:Case1c-g}, ~\eqref{trapezoid:Case1a-g}, 
and~\eqref{trapezoid:Case1b-g}.   

The major and minor halves of the 
center-justified Ferrers diagram for the partition~\eqref{trapezoid:UsualPartition}  of $n$ 
also have Durfee squares with side $j$.  
The associated major-minor hook partition for $n$, denoted 
\[
n = h_1 + h'_1 + \cdots + h'_{r_k-1} +  h_{r_k}  + h'_{r_k}  + h_{r_k + 1} + \cdots + h'_j 
\]
is a  partition into strictly decreasing parts, where 
\begin{align*}
h_i & = 
\begin{cases}
r_i + k - 2i + 1 & \text{for $i=1,\ldots,  r_k $} \\
  r_i + j - 2i + 1 & \text{for $i= r_k + 1 , \ldots, j$.}
  \end{cases}   \\
h'_i & = 
\begin{cases}
r_i + k -2i  & \text{for $i=1,\ldots, r_k - 1$} \\
  r_i + j -2i    & \text{for $i= r_k, \ldots, j$.}
\end{cases}
\end{align*}
Applying~\eqref{trapezoid:Case1a-g} 
and~\eqref{trapezoid:Case1b-g}, we obtain, for $i=1,\ldots, r_k - 1$,  
\[
h_i -  g_ i = h'_i -  g'_ i = h_{r_k} - g_{r_k} = k-j 
\]
and, for $i=r_k + 1,\ldots, j$,  
\[
h_i = g_ i = r_i + j - 2I+1 
\]
and
\[
h'_i = g'_ i  = r_i + j - 2I. 
\]
The critical observations are that 
\[
 h'_{r_k} =  g'_{r_k} = r_{r_k} + j -2r_k
\] 
\[
g_{r_k} - g'_{r_k} = (r_{r_k} + j - 2r_k + 1) - ( r_{r_k} + j -2r_k ) = 1
\]
and  
\begin{align*}
h_{r_k} - h'_{r_k} & =  \left(  r_{r_k} + k - 2 r_k + 1  \right)      
-  \left(  r_{r_k} + k - 2 r_k \right)     \\
& = k-j+1 \geq 2.
\end{align*}
These imply that the number of pairwise disjoint maximal sequences 
of consecutive integers in the major-minor hook partition for $n$ 
is exactly one more than the number in the  major-minor hook partition for $m$.
By the induction hypothesis,  the hook partition for $m$ contains exactly $\ell - 1$ 
such sequences, 
and so the hook partition for $n$ contains exactly $\ell$ 
pairwise disjoint maximal  sequences of consecutive integers.  

For example, if 
\[
n = 57 = 11 + 11 + 11 + 9 + 5 + 5 + 5
\]
then $k = 7$, $\ell = 3$ and 
\[
r_7 = 3 < j = 4 < r_j = 5.
\]
We have $k-j = 3$ and 
\[
m = 42  = 11 + 11 + 11 + 9.
\] 
The major-minor hook partition for $m$ is 
\[
m = 42 = 9 + 8 + 7 + 6 + 5 + 4  + 2 + 1 
\]
and contains two maximal sequences of consecutive integers:
$(9,8,7,6,5,4)$ and $(2,1)$.  
The major-minor hook partition for $n$ is   
\[
n = 57 = 12 + 11 + 10 + 9 + 8 + 4 + 2 + 1
\]
and contains three maximal sequences of consecutive integers:
$(12, 11, 10, 9, 8)$, $(4)$, and $(2,1)$.   
Note that $r_k = r_7 = 3$, $h_3 = 8$, $g_3 = 5$, and $h'_3 = g'_3 = 4$.

Case 3: 
\[
\boxed{r_k  < r_j = j}
\] 
Because $r_k = r_{j+1} < r_j = j$, it follows that the sides of the Durfee squares 
of the major halves of the partitions of both $m$ and $n$ are $s=j$.
Because $r'_j = r_j-1 = j-1$ and 
$r'_{j-1} = r_{j-1}-1 \geq r_j - 1 = j-1$, 
it follows that the sides of the Durfee squares 
of the minor halves of the partitions of both $m$ and $n$ are $s' = j-1$.
The hook numbers of the major halves are  
\begin{align*}
g_i & = r_i + j - 2i + 1  \qquad\text{for $i=1,\ldots, j$}\\
h_i & = \begin{cases}
r_i + k- 2i + 1 & \text{if $1 \leq i \leq r_k$}\\
r_i + j - 2i + 1 & \text{if $r_k+1 \leq i \leq j$}
\end{cases}
\end{align*}
and so 
\[
h_i - g_i = 
\begin{cases}
k-j  & \text{if $1 \leq i \leq r_k$}\\
0 & \text{if $r_k+1 \leq i \leq j$} 
\end{cases}
\]
The hook numbers of the minor halves are  
\begin{align*}
g'_i & = r_i + j - 2i   \qquad\text{for $i=1,\ldots, j-1$}\\
h'_i & = \begin{cases}
r_i + k- 2i  & \text{if $1 \leq i \leq r_k - 1$}\\
r_i + j - 2i & \text{if $r_k \leq i \leq j$}
\end{cases}
\end{align*}
and so 
\[
h'_i - g'_i = 
\begin{cases}
k-j  & \text{if $1 \leq i \leq r_k - 1$}\\
0 & \text{if $r_k \leq i \leq j$} 
\end{cases}
\]
Because 
\[
g_{r_k} - g'_{r_k} = (r_{r_k}+j-2r_k+1) - (r_{r_k}+j-2r_k) = 1
\]
and
\[
h_{r_k} - h'_{r_k} = (r_{r_k}+ k -2r_k+1) - (r_{r_k}+j-2r_k) = k-j + 1 \geq 2
\]
it follows that the major-minor hook partition for $n$ contains exactly more 
maximal sequence of consecutive integers than the hook partition for $m$.  

For example, if 
\[
n = 13 = 9 + 3 + 1
\]
then $k = \ell = 3$ and 
\[
r_3 = 1 < r_j = 2 = j.
\]
We have $k-j = 1$ and
\[
m = 12 = 9 + 3.
\]
The major-minor hook partition for $m$ is 
\[
m = 12 = 6 + 5 + 1
\]
and contains two maximal sequences of consecutive integers:
$(6,5)$ and $(1)$.  
The major-minor hook partition for $n$ is   
\[
n = 13 = 7 + 5  + 1
\]
and contains three maximal sequences of consecutive integers:
$(7)$, $(5)$, and $(1)$.

Case 4: 
\[
\boxed{r_k  < r_j < j}
\] 
Let $s$ be the side of the Durfee square of the major half of partition of $n$.  
The inequality $r_j < j$ implies that $s \leq j-1$, 
and so
\[
r_k < r_j \leq r_{s+1} \leq s < j.
\]
The side of the Durfee square of the major half 
of the partition of $m$ is also $s$.  
Let $c_i$ be the number of dots in the $i$th column of the Ferrers diagram of 
the major half of the partition of $n$.
We have
\[
g_i = 
\begin{cases}
r_i + j - 2i + 1 & \text{if $1 \leq i \leq r_j$}\\
r_i + c_i - 2i + 1 & \text{if $r_j +1 \leq i \leq s$}
\end{cases}
\]
and
\[
h_i = 
\begin{cases}
r_i + k - 2i + 1 & \text{if $1 \leq i \leq r_k  $}\\
r_i + j - 2i + 1 & \text{if $r_k + 1 \leq i \leq r_j$}\\
r_i + c_i - 2i + 1 & \text{if $r_j +1 \leq i \leq s$.}\\
\end{cases}
\]
Thus, 
\[
h_i - g_i = 
\begin{cases}
k - j & \text{if $1 \leq i \leq r_k$}\\
0 & \text{if $r_k +1 \leq i \leq s$.}\\
\end{cases}
\]

Let $c'_i = c_{i+1}$ be the number of dots in the $i$th column of the minor half 
of the partition of $n$
Let $s'$ denote the side of the Durfee square of the minor half of the partition of $n$.  
If $r_s \geq s+1$, then 
\[
r'_s = r_s - 1 \geq s \geq r_{s+1} > r'_{s+1}
\]
and so 
\[
s' = s.
\]
If $r_s = s$, then 
\[
r'_{s-1} = r_{s-1} - 1 \geq r_s - 1 = s-1 = r'_s
\]
and so 
\[
s' = s -1.
\]
In both cases we have $r_j - 1 \leq s-1 \leq s'$ and 
\[
g'_i = 
\begin{cases}
r_i + j - 2i & \text{if $1 \leq i \leq r_j - 1$}\\
r_i + c'_i - 2i & \text{if $r_j  \leq i \leq s'$}\\
\end{cases}
\]   
and
\[
h'_i = 
\begin{cases}
r_i + k - 2i  & \text{if $1 \leq i \leq r_k-1  $}\\
r_i + j - 2i  & \text{if $r_k  \leq i \leq r_j - 1$}\\
r_i + c'_i - 2i & \text{if $r_j +1 \leq i \leq s'$.}\\
\end{cases}
\]
Thus,  
\[
h'_i - g'_i  = 
\begin{cases}
k - j & \text{if $1 \leq i \leq r_k - 1$}\\
0 & \text{if $r_k  \leq i \leq s'$.}\\
\end{cases}
\]
Because 
\[
g_{r_k} - g'_{r_k} = (r_{r_k} +j - 2r_k + 1) - ( r_{r_k} + j -2r_k ) = 1
\]
and 
\[
h_{r_k} - h'_{r_k} =  (r_{r_k} + k  - 2r_k + 1) - ( r_{r_k} + j -2r_k ) 
= k-j+1 \geq 2
\]
it follows that the major-minor hook partition for $n$ contains exactly one more 
sequence of consecutive integers than the hook partition for $m$.  

For example, if 
\[
n = 50 = 11 + 11 + 9 + 7 + 7 + 5
\]
then $k=6$, $\ell = 4$, and 
\[
r_6 = 3 < r_5 = 4 < j = 5. 
\]
We have $k-j = 1$ and 
\[
m = 45 = 11 + 11 + 9 + 7 + 7.  
\]
The major-minor hook partition for $m$ is 
\[
m = 45 = 10 + 9 + 8 + 7 + 5 + 4 + 2
\]
and contains three pairwise disjoint maximal sequences of consecutive integers:
$(10, 9, 8, 7)$, $(5,4)$, and $(2)$.
The major-minor hook partition for $n$ is 
\[
n = 50  = 11 + 10 + 9 + 8  + 6 + 4 + 2
\]
and contains four disjoint maximal sequences:
$(11, 10, 9, 8)$, $(6)$, $(4)$, and $(2)$.

This completes the proof.

\end{proof}

\bt[Sylvester]          \label{trapezoid:theorem:10}
For all positive integers $n$ and $\ell$, 
\[
U_{\ell}(n) = V_{\ell}(n).
\]
\et

\begin{proof}
By Theorem~\ref{trapezoid:theorem:Sylvester-big},
Sylvester's one-to-one and onto function $f: \mcu(n) \rightarrow \mcv(n)$ 
maps $\mcu_{\ell}(n)$ into $\mcv_{\ell}(n)$.
We simply apply Lemma~\ref{trapezoid:lemma:counting} to complete the proof.  
\end{proof}

There are several recent proofs of Theorem~\ref{trapezoid:theorem:10}, for example, 
Andrews~\cite{andr66}, 
Andrews and Eriksson~\cite{andr-erik04}, and Hirschhorn~\cite{hirs74}.  
V. Ramamani and K. Venkatachaliengar~\cite{rama72} obtained a combinatorial proof.  
Their method is discussed in Andrews~\cite[pp. 448--449]{andr74} and~\cite[pp.  24--25] {andr98}.

For other recent work on trapezoidal numbers, see 
Apostol~\cite{apos03},  
Guy~\cite{guy82}, 
Leveque~\cite{leve50}, 
Moser~\cite{mose63}, Pong~\cite{pong07,pong09}, 
and Tsai and Zaharescu~\cite{tsai-zaha12b,tsai-zaha12a}.

\section{A problem}
An odd integer is an integer of the form $r + (r-1)$.  
Thus, a partition into odd parts is a partition into parts, 
each of which is a sum of two consecutive integers.  
A different generalization of Euler's theorem about partitions into odd parts 
would be a theorem about partitions into parts, each of which is a sum of $e$ 
consecutive integers, or, equivalently, a sum of $e$-trapezoids.   
Thus, we consider positive parts of the form 
\[
a_i = \sum_{j=0}^{e-1} ( r_i - j )
\]
with
\[
r_i \geq e
\]
and partitions of the form 
\[
n = \sum_{i=1}^k a_i 
= \sum_{i=1}^k \left(  \sum_{j=0}^{e-1} ( r_i -j ) \right). 
%= \sum_{j=0}^{s-1}   \left(  \sum_{i=1}^k ( r_i -j ) \right) 
\]
Interchanging summations, we obtain a partition of $n$ into $e$ parts, 
each of which inherits a well-defined partition:
\[
n = \sum_{j=0}^{e-1} n_j 
\]
where 
\[
n_j = \sum_{i=1}^k ( r_i -j ).
\]
Partitions into 2-trapezoids (that is, partitions intp odd numbers) are equinumerous 
with partitions into distinct parts.  What kind of partition are in one-to-one 
correspondence with partitions into $e$-trapezoids for $e \geq 3$?

\emph{Acknowledgement}.  I thank the referee for providing many references 
to the current literature on Sylvester's theorem.

\def\cprime{$'$} \def\cprime{$'$} \def\cprime{$'$}
\providecommand{\bysame}{\leavevmode\hbox to3em{\hrulefill}\thinspace}
\providecommand{\MR}{\relax\ifhmode\unskip\space\fi MR }
\providecommand{\MRhref}[2]{%
  \href{http://www.ams.org/mathscinet-getitem?mr=#1}{#2}
}
\providecommand{\href}[2]{#2}

\end{document}